\documentclass[a4paper]{article}
\usepackage[utf8]{inputenc}

\usepackage[top=2.5cm, bottom=2.5cm, left=2.5cm, right=2.5cm]{geometry}

\usepackage{amsmath, amssymb, url, mathtools}
\usepackage{mathrsfs}
\usepackage[hidelinks]{hyperref}
\usepackage{tikz}

\usepackage[shortlabels]{enumitem}

 \usepackage{xcolor}

 \usepackage{amsthm}
 \usepackage[noabbrev]{cleveref}
 \newtheorem{thm}{Theorem}[section]
 \newtheorem{lm}[thm]{Lemma}
 \newtheorem{res}[thm]{Result}
 
 \newtheorem{prop}[thm]{Proposition}

 \theoremstyle{definition}

 \newtheorem{df}[thm]{Definition}

 \newcommand{\eps}{\varepsilon}

 \newcommand{\FF}{\mathbb F}

 \newcommand{\vspan}[1]{\left \langle #1 \right \rangle}
 \newcommand{\set}[1]{ \left \{ #1 \right \} }
 \newcommand{\sett}[2]{ \left\{ #1 \, \, || \, \, #2 \right \} }

 \newcommand{\pg}{\textnormal{PG}}
 \newcommand{\PG}{\textnormal{PG}}
 \newcommand{\one}{\mathbf 1}
 \newcommand{\zero}{\mathbf 0}
 \newcommand{\mq}{
 Q}
 
 \newcommand{\mc}{\mathcal C}
 \newcommand{\mg}{\mathcal G}
 \newcommand{\mr}{\mathcal R}

 \DeclareMathOperator{\supp}{supp}
 \DeclareMathOperator{\wt}{wt}
 \DeclareMathOperator{\proj}{proj}

\newcommand{\old}[1]{}

\title{The minimum weight of the code of intersecting  lines in $\pg(3,q)$}
\author{
 Sam Adriaensen \\ \textit{Vrije Universiteit Brussel} \and
 Robin Simoens \\ \textit{Ghent University} \\ \textit{Universitat Politècnica de Catalunya} \\ \and 
 Leo Storme \\ \textit{Ghent University}
}
\date{}

\begin{document}

\maketitle

\begin{abstract}
 We characterise the minimum weight codewords of the $p$-ary linear code of intersecting lines in $\pg(3,q)$, $q=p^h$, $q\geq19$, $p$ prime, $h\geq 1$.
 If $q$ is even, the minimum weight equals $q^3+q^2+q+1$. If $q$ is odd, the minimum weight equals $q^3+2q^2+q+1$.
 For $q$ even, we also characterise the codewords of second smallest weight.
\end{abstract}

\paragraph{Keywords.} Linear codes; Finite projective spaces; Klein quadric.

\paragraph{MSC.}
05B25, 
51A50, 
51E21, 
94B05. 

\section{Introduction}

Coding theory is the study of reliable communication through a noisy channel.
The challenge of coding theory is to construct codes with little redundancy that have strong error-correcting capabilities.
Finite geometry has contributed in this endeavour in several ways, one of which is explored further in this article.
More specifically, we investigate linear codes constructed from certain matrices arising from finite geometries.
Almost 60 years ago, Weldon \cite{weldon1966} initiated the study of such codes by analysing codes constructed from the incidence matrix of classical projective planes\footnote{Actually, Weldon \cite{weldon1966} studied codes from projective planes arising from perfect difference sets. Classical projective planes are known to belong to this category. It is still an open conjecture that no other projective planes do. This conjecture has been verified for planes up to order 2 billion \cite{BaumertGordon}.}, because they allow an efficient decoding algorithm.
More general classes of codes arising from finite geometries have been investigated since, see e.g.\ \cite{LavrauwC} for a detailed survey.
Most attention has been paid to codes arising from the incidence matrix of a finite projective space.
This works as follows.
Choose integers $0 \leq j < k < n$ and a prime power $q$.
Construct a $01$-matrix $M$ whose rows and columns are indexed by the subspaces of the projective space $\pg(n,q)$ of dimension $k$ and $j$ respectively.
Let $M$ have entry 1 in a certain position if and only if the corresponding $j$- and $k$-space are incident.
Then the code to be studied is the row space of $M$ over a finite field, typically the prime subfield of $\FF_q$.
The minimum weight of this code was determined by Delsarte, Goethals, and MacWilliams \cite{delsarte} in case $j=0$ and by Bagchi and Inamdar \cite{bagchi} for general values of $j$.
Stronger characterisations of the codewords of small weight are known, see e.g.\ \cite{adriaensen2023small}.
The dimension of these codes was determined in case $j=0$ by Hamada \cite{hamada}, but is still unknown if $j>0$ and $k < n-1$.

In this paper, we study a similar code.
Choose integers $0 \leq j,k < n$ and a prime power $q$.
Construct a $01$-matrix $M$ whose rows and columns are indexed by the subspaces of $\pg(n,q)$ of dimensions $k$ and $j$ respectively.
Let $M$ have entry 1 in a certain position if and only if the corresponding $j$- and $k$-space have a non-empty intersection.
Let $\FF_p$ denote the prime subfield of $\FF_q$ and $\mc_{j,k}(n,q)$ the row space of $M$ over $\FF_p$.
The length of this code equals the number of $j$-spaces of $\pg(n,q)$.
The dimension of this code was determined by Sin \cite{Sin:04}.
The general formula is quite involved, so we will not state it here.
However, the case of interest to us are the codes $\mc_{1,1}(3,q)$, and for this fixed set of dimensions $j,k,n$, there is a simple formula for the dimension of the code.

\begin{res}[{\cite{Sin:04}}]
 \label{Res:DimCode}
 Let $p$ be a prime and $q = p^h$ with $h \geq 1$ an integer.
 Then
 \[
  \dim \mc_{1,1}(3,q) = q \left( \frac{2p^2+1}3 \right)^h + 1.
 \]
\end{res}

Yet, determining the minimum weight of $\mc_{j,k}(n,q)$ is still an open problem for $j,k>0$ and $j+k < n-1$ \cite[Open Problem 2.4]{LavrauwC}.
In this paper, we determine the minimum weight of $\mc_{1,1}(3,q)$ for sufficiently large $q$.
For a set $S$ of lines in $\pg(3,q)$, we call the $01$-vector whose positions are indexed by the lines of $\pg(3,q)$ and that has a 1 exactly in the positions corresponding to the lines of $S$, the characteristic vector of $S$.
For any line $\ell$ of $\pg(3,q)$, we refer to the characteristic vector of the set of lines having non-empty intersection with $\ell$ simply as the characteristic vector of $\ell$.
Our main result is the following.

\begin{thm}
 \label{Thm:Main:Code}
 Suppose that $q \geq 19$ and let $c$ be a non-zero codeword of $\mc_{1,1}(3,q)$ of weight at most $q^3+2q^2+q+1$.
 Then either
 \begin{enumerate}
  \item \label{Thm:Main:Code:Symplectic} $c$ has weight $q^3+q^2+q+1$ and $c$ is the characteristic vector of the set of absolute lines of a symplectic polar space $W(3,q)$,
  \item $c$ has weight $q^3+2q^2+q+1$ and $c$ is a scalar multiple of the characteristic vector of a line.
 \end{enumerate}
 Moreover, (\ref{Thm:Main:Code:Symplectic}) occurs if and only if $q$ is even.
\end{thm}

In order to prove this theorem, we use a classification of small blocking sets on the Klein quadric.
Metsch \cite{metsch2000blocking} proved the following result, where $\delta$ denotes the Kronecker delta.

\begin{res}[{\cite[Theorem 1.1]{metsch2000blocking}}]
 \label{Res:Metsch}
 Let $B$ be a minimal blocking set with respect to lines on the hyperbolic quadric $Q^+(n,q)$, with $n \geq 5$.
 If
 \[
  |B| \leq \frac{q^{n-1}-1}{q-1} + q^\frac{n-1}2 - \delta_{n,5},
 \]
 then $B$ is contained in a hyperplane section of $Q^+(n,q)$.
\end{res}

It is easy to check that a minimal blocking set $B$ with respect to the lines of $Q^+(n,q)$, which is contained in a hyperplane section, can take two forms:
 \begin{enumerate}
  \item $B$ is a parabolic hyperplane section of $Q^+(n,q)$, in which case $|B| = \frac{q^{n-1}-1}{q-1}$.
  \item There exists a point $P \in Q^+(n,q)$ with tangent hyperplane $P^\perp$ such that $B$ consists of the hyperplane section $P^\perp \cap Q^+(n,q)$, where the point $P$ is removed.
  In this case $|B| = \frac{q^{n-1}-1}{q-1} + q^\frac{n-1}2 - 1$.
 \end{enumerate}

We are specifically interested in the case $n=5$, but we need to strengthen the above result slightly by removing the term $\delta_{n,5}$ from the upper bound on $|B|$.
We will do this on the condition that $q \geq 4$, and prove the following theorem.

\begin{thm}
 \label{Thm:Main:Klein}
 Suppose that $q \geq 4$.
 Let $B$ be a minimal blocking set with respect to the lines of $Q^+(5,q)$.
 If $|B| \leq q^3+2q^2+q+1$, then $B$ is contained in a hyperplane section of $Q^+(5,q)$.
\end{thm}

\paragraph{Overview.} In \Cref{Sec:Prel}, we state the necessary background.
In \Cref{Sec:Klein}, we prove \Cref{Thm:Main:Klein}.
Finally, in \Cref{Sec:MinWt}, we prove \Cref{Thm:Main:Code}.

\section{Preliminaries}
 \label{Sec:Prel}

Throughout this article, $p$ denotes a prime number, $h\geq1$ an integer, and $q = p^h$ a prime power.
The finite field of order $q$ will be denoted as $\FF_q$.

\bigskip

We recall the basics of the theory of linear codes, see e.g.\ \cite[\S 3]{vanLint99}.
It will be convenient to represent the ambient vector space of a linear code as the vector space $\FF_q^S$ of functions from some finite set $S$ to $\FF_q$.
Given a vector $v \in \FF_q^S$, its \emph{support} is defined as
\[
 \supp(v) = \sett{ s \in S }{ v(s) \neq 0 },
\]
and its \emph{(Hamming) weight} as $\wt(v) = |\supp(v)|$.
The \emph{(Hamming) distance} between two vectors $v$ and $w$ is $d_H(v,w) = \wt(v-w)$.
A subspace $C$ of $\FF_q^S$ is said to be a \emph{linear $[n,k,d]_q$ code} if
\begin{itemize}
 \item its \emph{length}, which equals $|S|$, is $n$,
 \item its dimension is $k$,
 \item its \emph{minimum distance} $\min \sett{d_H(v,w)}{ v,w \in C, \, v \neq w}$ equals $d$.
\end{itemize}
It is well-known that the minimum distance of $C$ equals its \emph{minimum weight}, which is defined as $\min \sett{ \wt(v) }{ v \in C \setminus \{\zero\} }$, where $\zero$ denotes the zero vector.

The standard scalar product on $\FF_q^S$ is defined as follows, for $v, w \in \FF_q^S$:
\[
 v \cdot w = \sum_{s \in S} v(s) w(s).
\]
If $C \subseteq \FF_q^S$ is a linear code, then its orthogonal complement with respect to the standard scalar product is denoted as $C^\perp$ and is called the \emph{dual code} of $C$.

\bigskip

The projective space corresponding to the vector space $\FF_q^{n+1}$ will be denoted as $\pg(n,q)$.
From now on, all dimensions will be projective dimensions. A space of dimension $k$ is called a $k$-space for short.
Let $\mg_k(n,q)$ denote the set of $k$-spaces of $\pg(n,q)$.
For a subset $S$ of $\mg_k(n,q)$, define its \emph{characteristic vector} as
\[
 \chi_S: \mg_k(n,q) \to \{0,1\}: \kappa \mapsto \begin{cases}
  1 & \text{if } \kappa \in S, \\
  0 & \text{if } \kappa \notin S.
 \end{cases}
\]
For an integer $j$ and a $k$-space $\kappa$ of $\pg(n,q)$, let $\chi_\kappa^{(j)}$ denote the characteristic vector of the set of $j$-spaces of $\pg(n,q)$ having non-empty intersection with $\kappa$.
We refer to $\chi_\kappa^{(j)}$ as the characteristic vector of $\kappa$.

\begin{df}
 Let $q$ be a power of a prime $p$, and choose integers $0 < j,k < n$.
 The code $\mc_{j,k}(n,q)$ is defined as the $\FF_p$-span of $\sett{\chi_\kappa^{(j)}}{\kappa \in \mg_k(n,q)}$.
\end{df}

The code of interest to us is $\mc_{1,1}(3,q)$, which we will simply denote as $\mc(3,q)$.
The code $\mc_{0,1}(2,q)$ also plays an important role in this article, as we will heavily rely on the characterisation of its small weight codewords by Sz\H onyi and Weiner \cite{szonyiB}.

\begin{res}[{\cite[\S 4]{szonyiB}}]
 \label{Res:PlaneCodes}
 Suppose that $q \geq 19$ and that $c \in \mc_{0,1}(2,q)$ with
 \[
  \wt(c) < \begin{cases}
   3q-3 & \text{if $q$ is prime}, \\
   3q-12 & \text{otherwise.}
  \end{cases}
 \]
 Then $c$ is a linear combination of at most 2 characteristic vectors of lines.
\end{res}

We end the preliminaries by discussing quadrics and symplectic polar spaces.
We refer the reader to \cite[\S 1]{hirschfeldthas} for an in-depth treaty of quadrics, and to \cite{Cameron15} for a more general treatment of polar spaces.
A form $b: \FF_q^{n+1} \times \FF_q^{n+1} \to \FF_q$ is called \emph{bilinear} if it is linear in both arguments.
We call a bilinear form $b$ \emph{reflexive} if either
\begin{itemize}
 \item $b$ is \emph{symmetric}, i.e.\ $b(x,y) = b(y,x)$ for all $x,y \in \FF_q^{n+1}$, or
 \item $b$ is \emph{alternating}, i.e.\ $b(x,y) = -b(y,x)$ and $b(x,x) = 0$ for all $x,y \in \FF_q^{n+1}$.
\end{itemize}
Now suppose that $b$ is a reflexive bilinear form. 
A non-zero vector $x \in \FF_q^{n+1}$ is called \emph{singular} with respect to $b$ if $b(x,y) = 0$ for all $y\in \FF_q^{n+1}$.
We call $b$ \emph{degenerate} if it has a singular vector.
If $b$ is non-degenerate, then it induces an involution $\perp$ on the subspaces of $\FF_q^{n+1}$ defined by
\[
 \pi^\perp = \sett{y \in \FF_q^{n+1}}{(\forall x \in\pi)(b(x,y) = 0)}.
\]
Since the subspaces of $\pg(n,q)$ are essentially the subspaces of $\FF_q^{n+1}$, we can interpret $\perp$ as being defined on $\pg(n,q)$.
Note that $\perp$ reverses inclusion.
Any involution on the subspaces of $\pg(n,q)$ that reverses inclusion is called a \emph{polarity}.
A subspace $\pi$ of $\pg(n,q)$ is called \emph{absolute} with respect to the polarity $\perp$ if $\pi \subseteq \pi^\perp$.

If $b$ is a non-degenerate alternating form on $\FF_q^{n+1}$, then $n$ is necessarily odd, and up to a change of basis, $b$ is given by
\[
 b(x,y) = x_0 y_1 - x_1 y_0 + \cdots + x_{n-1} y_n - x_n y_{n-1}.
\]
We call the corresponding polarity $\perp$ \emph{symplectic}, and the set of absolute subspaces with respect to $\perp$ is called the \emph{symplectic polar space} $W(n,q)$.
It is worth noting that every point of $\pg(n,q)$ is absolute with respect to a given symplectic polarity $\perp$, and that a line $\ell$ through a point $P$ is absolute if and only if $\ell \subset P^\perp$.

A form $f: \FF_q^{n+1} \to \FF_q$ is called \emph{quadratic} if $f(\alpha x) = \alpha^2 f(x)$ for all $x \in \FF_q^{n+1}$ and all $\alpha \in \FF_q$.
A non-zero vector $x$ is called \emph{isotropic} with respect to $f$ if $f(x)=0$ and \emph{anisotropic} otherwise.
A point $P \in \pg(n,q)$ is called (an)isotropic if its coordinate vectors are (an)isotropic.
The set of isotropic points of $\pg(n,q)$ with respect to a quadratic form is called a \emph{quadric}.
Every quadratic form $f$ gives rise to a symmetric bilinear form $b(x,y) = f(x+y) - f(x) - f(y)$.
We call $f$ and its corresponding quadric \emph{degenerate} if some isotropic vector $x$ is singular with respect to $b$.
The non-degenerate quadratic forms on $\FF_q^{n+1}$ have been classified.
\begin{itemize}
 \item If $n$ is odd, there are, up to a change of basis, two non-degenerate quadratic forms on $\FF_q^{n+1}$.
 They are of the form
 \[
  f(x) = x_0 x_1 + \cdots + x_{n-3} x_{n-2} + g(x_{n-1},x_n)
 \]
 for some non-degenerate quadratic form $g$.
 The corresponding quadric is the \emph{hyperbolic quadric} $Q^+(n,q)$ or the \emph{elliptic quadric} $Q^-(n,q)$, if $g$ is reducible or irreducible over $\FF_q$ respectively.
 The corresponding bilinear form $b$ is non-degenerate and defines a polarity $\perp$.
 If $q$ is odd, then the set of absolute points with respect to $\perp$ coincides with the quadric $\mq^\pm(n,q)$.
 If $q$ is even, then $\perp$ is a symplectic polarity.
 \item If $n$ is even, there is up to a change of basis a unique non-degenerate quadratic form on $\FF_q^{n+1}$, given by
 \[
  f(x) = x_0 x_1 + \cdots + x_{n-2} x_{n-1} + x_n^2.
 \]
 The corresponding quadric is called \emph{parabolic} and denoted as $\mq(n,q)$.
 Let $b$ denote the corresponding bilinear form.
 Then $b$ is non-degenerate if and only if $q$ is odd, in which case $\mq(n,q)$ is the set of absolute points of the corresponding polarity.
\end{itemize}
We also use the notation $Q^\eps(n,q)$, where $\eps$ is either $-1$, $0$, or $+1$, in which case the quadric is respectively elliptic, parabolic, or hyperbolic. In other words, $Q^{-1}(n,q)=Q^-(n,q)$, $Q^0(n,q)=Q(n,q)$ and $Q^1(n,q)=Q^+(n,q)$.

Let $\pi$ and $\rho$ be two disjoint subspaces of $\pg(n,q)$ and let $S$ be a set of points in $\pi$.
The set $\rho S = \bigcup_{P \in \pi} \vspan{\rho,P}$ is called the \emph{cone} with \emph{base} $S$ and \emph{vertex} $\rho$.
By convention, $\rho S = S$ if $\rho = \emptyset$ and $\rho S = \rho$ if $S = \emptyset$.
Every quadric is a cone with as base a non-degenerate quadric.
We call a subspace $\pi$ \emph{singular} with respect to a quadric $Q$ if $\pi \cap Q$ is a degenerate quadric.
If $P \in Q^\eps(n,q)$ and $Q^\eps(n,q)$ has a polarity $\perp$, then $P^\perp$ intersects $Q^\eps(n,q)$ in a quadric $P Q^\eps(n-2,q)$.
We call $P^\perp$ the \emph{tangent hyperplane} of $P$.
A subspace contained in $Q^\eps(n,q)$ of maximal dimension is called a \emph{generator}.
\begin{res}
 \label{Res:Quadrics}
 \begin{enumerate}
  \item The quadric $Q^\eps(n,q)$ has $\frac{q^n-1}{q-1} + \eps q^\frac{n-1}2$ points.
  \item \label{Res:Quadrics:Generators} The quadric $Q^\eps(n,q)$ has
  \[
   \prod_{i=1-\eps}^\frac{n-\eps}2 (q^i+1)
  \]
  generators, and their dimension equals $\frac{n+\eps}2-1$.
  \item \label{Res:Quadrics:Equivalence} Let $\mg$ denote the set of generators of $Q^+(n,q)$.
  Then the relation
  \[
   \sett{ (\pi, \rho) \in \mg^2 }{ \dim (\pi \cap \rho) \equiv \frac{n-1}2 \pmod 2 }
  \]
  is an equivalence relation on $\mg$, having two classes of equal size.
 \end{enumerate}
\end{res}

An equivalence class of generators of $Q^+(3,q)$ is called a \emph{regulus}, and two distinct reguli of the same quadric $Q^+(3,q)$ are called \emph{opposite}.

Lastly, we discuss the \emph{Klein correspondence}.
This is a bijection between the points of $Q^+(5,q)$ and the lines of $\pg(3,q)$.
For this reason, the quadric $Q^+(5,q)$ is also called the \emph{Klein quadric}.
The lines contained in $Q^+(5,q)$ correspond to the sets $\sett{\ell \in \mg_1(3,q)}{P \in \ell \subset \pi}$ of lines in $\pg(3,q)$, with $P$ and $\pi$ an incident point and plane of $\pg(3,q)$.
The equivalence classes of generators of $Q^+(5,q)$ are called the classes of \emph{Greek} and \emph{Latin} planes.
The Greek planes correspond to the sets $\sett{\ell \in \mg_1(3,q)}{\ell \subset \pi}$ with $\pi$ a plane, and the Latin planes correspond to the sets $\sett{\ell \in \mg_1(3,q)}{P \in \ell}$ with $P$ a point.
A parabolic hyperplane section of $Q^+(5,q)$ corresponds to the set of lines of a symplectic polar space $W(3,q)$.
A tangent hyperplane section of $Q^+(5,q)$ corresponds to the set of all lines having non-empty intersection with a fixed line of $\pg(3,q)$.

\section{Blocking sets of the Klein quadric}
 \label{Sec:Klein}

In this section, we extend \Cref{Res:Metsch} in the case $n=5$ to \Cref{Thm:Main:Klein}.
We will closely follow the original proof by Metsch \cite{metsch2000blocking}.

\begin{df}
 A \emph{blocking set} of a quadric $Q$ with respect to the $k$-spaces is a set $B$ of points of $Q$ that intersects every $k$-space of $Q$.
 We call $B$ \emph{minimal} if none of its proper subsets is a blocking set.
\end{df}

From now on, blocking sets are always considered with respect to lines.
We focus on $Q^+(5,q)$, which will simply be denoted as $Q$.
The corresponding polarity will be denoted as $\perp$.

\bigskip

Metsch proved the following result.

\begin{res}[{\cite[Lemma 2.2]{metsch2000blocking}}]
 \label{Res:Metsch:Lemma2.2}
 Let $B$ be a blocking set of $Q$ with $|B| \leq q^3+2q^2+q$.
 Let $H$ be a hyperplane of $\pg(5,q)$ such that $H \cap B$ is not a blocking set of $Q$.
 \begin{enumerate}
  \item If $H$ is a tangent hyperplane, then $|H \cap B| \leq |B| - q(q^2-q+1) \leq 3q^2$.
  \item If $H$ is a non-tangent hyperplane, then $|H \cap B| \leq 4q^2-2q$.
  If in addition $|H \cap B| > 3q^2-q$, then every line of $H \cap Q$ meets $B$ in at least two points.
 \end{enumerate}
\end{res}

We extend this result using similar arguments to the following lemma.

\begin{lm}\label{lemma:block}
 Suppose that \(q \geq 3\).
 Let \(B\) be a blocking set of \(Q\) with \(|B|\leq q^3+2q^2+q+1\).
 Let \(H\) be a hyperplane of \(\PG(5,q)\) such that \(H\cap B\) is not a blocking set of \(Q\).
 If \(|H\cap B|\geq 3q^2+2\), then \(H\) is a non-tangent hyperplane, \(|H \cap B| \leq 4q^2-2q+1\), and every line of \(H\cap Q\) meets \(B\) in at least two points.
\end{lm}

\begin{proof}
Note that the case \(|B| \leq q^3+2q^2+q\) follows from \Cref{lemma:block}, so we may suppose that $|B| = q^3+2q^2+q+1$.

In case that $H$ is a non-tangent hyperplane, put ${\mathcal P} \coloneqq H \cap Q$. In case that $H$ is the
tangent hyperplane of a point $P_0 \in Q$, put ${\mathcal P} \coloneqq (H \cap Q) \setminus  \{P_0\}$. In both cases, ${\mathcal P}$  is a minimal blocking set of $Q$.

There exist integers $\delta, c_0, c_1, d, \Delta$ such that the following properties hold:
\begin{enumerate}[(a)]
    \item  every point $P \in {\mathcal P}$ belongs to $c_0$ lines of $Q$ that do not lie in $H$;
    \item every line contained in ${\mathcal P}$ lies in $c_1$ planes of $Q$ that do not lie in $H$;
    \item $\delta$ points of ${\mathcal P}$ are not in $B$;
    \item $d$ is the maximal number of points of $Q \setminus B$ on a line contained in ${\mathcal P}$;
    \item $\Delta= |B| - |{\mathcal P}|$, which implies that $|B \setminus H| = |B| - |B \cap H | \leq  \delta + \Delta$.
\end{enumerate}

Note that $|B\setminus H| = \delta + \Delta$ or $|B \setminus H| = \delta + \Delta - 1$ and that the latter occurs only when $H$
is the tangent hyperplane of a point $P_0 \in  B$.

Since $B$ meets every line of $Q$, $d\neq  q + 1$. The definition of ${\mathcal P}$ implies that ${\mathcal P}$ meets
every line of $Q$.
Since, by hypothesis, the points of $H \cap  B$ do not meet every line of $Q$, we see that ${\mathcal P}$ is not contained in $B$. Hence, $d \neq  0$. Thus $1 \leq  d \leq  q$.

If $X$ is a point of $B$ that does not belong to $H$, then $X^\perp$ meets $H$ in a quadric $Q_X$ of type $Q^+(3,q)$.
Denote by $z$ the maximal number of points of $Q \setminus B$ that lie in such a quadric $Q_X$.  
In case
that $H$ is the tangent hyperplane of a point $P_0 \in Q$, note that $P_0$ does not occur in any of
the quadrics $Q_X$.

Then every point $X$ of $B \setminus H$ is perpendicular to at most $z$ points of ${\mathcal P} \setminus  B$. On the other
hand, a point $P$ of ${\mathcal P} \setminus  B$ belongs to $c_0$ lines of $Q$ that are not contained in $H$ and as each of
these lines meets $B$, the point $P$ is perpendicular to at least $c_0$ points of $B \setminus  H$. 
We now count the pairs $(X, P)$ of perpendicular points $X \in  B \setminus  H$ and $P \in {\mathcal P}  \setminus  B$. We obtain that $|B \setminus  H|\cdot z \geq \delta c_0$.

Since $|B \setminus H | \leq \Delta+\delta$, this implies that

\begin{equation}
(\Delta+\delta)\cdot z \geq \delta c_0.
\end{equation}

Consider a point $X \in B \setminus H$ for which the quadric $Q_X$ has exactly $z$ points in $Q \setminus  B$. Since
every line of $Q_X$ meets $Q \setminus B$ in at most $d$ points, a simple double counting argument (cf.\ \cite[Lemma~2.1]{metsch2000blocking}) gives
\[z \leq \frac{|Q_X|d}{q+1}.\]

So this   upper and lower bound for $z$ imply that
\begin{equation}\label{eq:2.2}
\frac{\delta c_0}{\delta+\Delta} \leq z \leq \frac{|Q_X|d}{q+1}.
\end{equation}

By definition (d) of $d$, ${\mathcal P}$ contains a line $\ell$ with exactly $d$ points $P_1,\dots,  P_d$ in $Q \setminus B$. The point $P_1$
belongs to $c_0$ lines of $Q$ that do not lie in $H$ and each of these lines meets $B \setminus  H$. In this way, the
point $P_1$ gives rise to at least $c_0$ points of $B\setminus H$. For the point $P_2$, the same holds, but $c_1q$ of these
lines on $P_2$ lie in a plane of $Q$ on $Q^+(5,q)$ through $\ell$, so that the corresponding points of $B \setminus H$ might have 
already been counted using the lines on $P_1$. Thus $P_2$ gives rise to at least $c_0 - c_1q$ new points in
$B \setminus  H$ and the same is true for $P_3,\dots, P_d$.  This shows that

\begin{equation}\label{eq:2.3}
|B \setminus  H| \geq dc_0 - (d - 1)c_1 q.
\end{equation}

 But $|B \setminus  H| \leq  \Delta + \delta$, so
$\delta \geq  dc_0 - (d - 1)c_1q - \Delta$.

Inequality (\ref{eq:2.2}) remains true, if we replace $\delta$ by this lower bound for $\delta$. This gives, after
simplifications using that $|Q_X|=(q+1)^2$, that
\begin{equation}\label{eq:2.4}
(dc_0- (d - 1)c_1q-\Delta) c_0\leq  (q + 1)  d(dc_0-(d-1)c_1q).
\end{equation}

There are two possibilities for the intersection of the hyperplane $H$ with the quadric $Q$.
Either $H$ intersects $Q$ in a non-degenerate parabolic quadric $Q(4,q)$ or $H$ is the tangent hyperplane $P^\perp$ of a point $P \in Q$.
    
Suppose by contradiction that $H$ is a tangent hyperplane. Then $c_0=q^2$ and $c_1=1$. The bound on $|B|$ gives $\Delta\leq 1$. From (\ref{eq:2.4}), we obtain $ \left(dq^2-(d-1)q\right)\left(q^2(q+1)-d(q+1)^2\right)\leq q^2(q+1)$, so $d=q$. Then $|B \setminus  H| \geq dc_0 - (d - 1)c_1 q=q^3-q^2+q$,  and thus $|B\cap H|\leq 3q^2+1$, a contradiction.

Therefore, $H$ is a non-tangent hyperplane. Then $c_0=q^2+q$ and $c_1=2$. The bound on $ |B|$ gives $\Delta\leq q^2$. From (\ref{eq:2.4}), we obtain $\left(dq-d+2\right)\left(q-d\right)\leq q^2$, so either \(d=1\), \(d=q-1\) or \(d=q\).

If $d=1$, then (\ref{eq:2.2}) implies that $\delta\leq q+1$ and therefore $ |B\setminus H|\leq q^2+q+1$. Also, (\ref{eq:2.3}) implies that $ |B\setminus H|\geq q^2+q$ and therefore $\delta\geq q$. The latter means that there are at least $q$ points in $H \cap Q$  not in $B$, and since $d=1$, these points are pairwise non-collinear in $Q$. Take two of these points. The number of points of $B\setminus H$ that are collinear with both points is at least $q^2+q-1$, since each point belongs to $c_0=q^2+q$ lines of $Q^+(5,q)\setminus H$, and $ |B\setminus H|\leq q^2+q+1$. There are at most $2$ points not collinear with both. The $q^2+q-1$ points that are, belong to a hyperbolic quadric $Q^+(3,q)$. If we take a third point of $B\setminus H$, then it can be collinear with at most $2q+1$ points of that hyperbolic quadric $Q^+(3,q)$, i.e.\ at most $2q+1+4=2q+5$ lines of the $c_0=q^2+q$ lines of $Q^+(5,q)\setminus H$ through that third point can be blocked by $B$, a contradiction.

If $d=q-1$, then (\ref{eq:2.3}) implies that $ |B\setminus H|\geq q^3-2q^2+3q$, so $ |B\cap H|\leq 4q^2-2q+1$. Moreover, $d=q-1$ means that every line of $H\cap Q$ contains at least two points of $B$.

If $d=q$, then (\ref{eq:2.3}) implies that $|B\setminus H|\geq q^3-q^2+2q$, so $|B\cap H|\leq 3q^2-q+1$, a contradiction.
\end{proof}

\begin{res}[{\cite[Lemma 2.4]{metsch2000blocking}}]
 \label{lemma:solid}
Consider in $\PG(4,q)$ a degenerate quadric $Q' = P Q^+(3,q)$ for some point $P$.
Assume that $B$ is a set of at most $ (q+1)^2+(2q-3)$ points of $Q'$ such that $P\notin B$ and such that every line of $Q'$ meets $B$. Then there exists a solid $S$ with $P\notin S$ and
$S\cap Q'\subseteq B$. 
\end{res}

We are now ready to prove \Cref{Thm:Main:Klein}, following a similar strategy as \cite[Theorem 3.1]{metsch2000blocking}.

\begin{proof}[Proof of \Cref{Thm:Main:Klein}]
If $ |B|\leq q^3+2q^2+q$, then this is follows from \Cref{Res:Metsch}, so suppose that $ |B|=q^3+2q^2+q+1$. 
    
Let ${\mathcal F}$ be one of the two equivalence classes of planes of $Q$. Then $ {|\mathcal F}| = q^3 + q^2 + q + 1$, and
every point of $Q$ belongs to $q + 1$ planes of ${\mathcal F}$. It follows that
$\sum_{\pi \in {\mathcal F}} |\pi \cap B| = |B|(q + 1)$. Since
$|B| \leq  q(q + 1)^2+1$, this implies that there exists a plane $\pi\in {\mathcal F}$ with $|\pi \cap B|< q + 2$. It follows
 from a result of Bruen \cite{bruen} that $|\pi \cap B| = q + 1$ and that $\ell \coloneqq \pi\cap B$ is a line.
 
Consider a point $P$ of $Q \setminus  B$. Then \Cref{lemma:solid} shows that $|P^\perp \cap  B| \geq  (q + 1)^2$. It also shows
that $|P^\perp \cap  B | < (q + 1)^2 + 2q - 2$ implies that there exists a solid $S$ in the hyperplane
$P^\perp$, with $P \not\in S$ and $S \cap Q\subseteq B$. Since $P \not\in S$, the subspace $S$ meets $Q$ in a quadric of type
$Q^+(3, q)$, which has $ (q+ 1)^2$ points. For every point $P \in Q \setminus B$, put $e_P \coloneqq |P^\perp \cap   B | - (q+ 1)^2$.

Now we prove the following results.

\bigskip

(1) The set $\pi\setminus \ell$ contains at least one point $P$ with $e_P \leq   q - 1$, at least two points $P$ with
$e_P \leq  q$, and, if $q \geq 4$, at least $q + 2$ points $P$ with $e_P < 2q - 2$, unless $e_P=q$ for all the points $P\in\pi\setminus \ell$.

\bigskip

Recall that $\pi\cap B = \ell$. Every point $P$ of $\pi\setminus \ell$ is perpendicular to $|P^\perp\cap B|=e_P+(q+1)^2$ points of $B$. But then $P$ is collinear to $e_P+(q+1)^2-(q+1)=q^2 + q + e_P$ points of $B \setminus \ell$. Then every point of $B \setminus  \ell$  is perpendicular to $q$ points of $\pi\setminus \ell$, unless this point of $B$ lies in the second plane of $Q$ through the line $\ell$.

Therefore

\[\sum_{P\in\pi\setminus \ell}(q^2+q+e_P)\leq (|B|-q-1)q.
\]
So
\[\sum_{P\in \pi\setminus \ell}e_P\leq (|B|-q-1)q-(q^2+q)q^2\leq q^3,\]
since $|B|\leq q^3+2q^2+q+1$.

Since $|\pi\setminus \ell|= q^2$ and since the upper bound is $q^3$, there exists at least one point $P\in\pi\setminus \ell$ with $e_P\leq q-1$, unless $e_P=q$ for every point $P\in\pi\setminus \ell$. It is also valid that there exist at least two points $P\in\pi\setminus \ell$ with
$e_P \leq  q$, and, if $q \geq 4$, at least $q + 2$ points $P\in\pi\setminus \ell$ with $e_P < 2q - 2$.

\bigskip

   {\bf Case 1}. Assume that it is not so that for all  points $P\in\pi\setminus \ell$, $e_P=q$. \bigskip

   By (1), we find $q + 2$ points $P_1,\dots, P_{q+2}$  in $\pi \setminus  \ell$,  with $e_{P_i} < 2q - 2$. For each of these points
$P_i$, the hyperplane $P_i^\perp$ contains a solid $S_i$, with $P_i \not\in S_i$, and such that all the points of the quadric
$S_i \cap Q$ belong to  $B$. These quadrics $S_i \cap  Q$ are hyperbolic quadrics  $Q^+(3, q)_i$.

Here, $\pi\subseteq P_i^\perp$. Since all the intersection points of the plane $\pi$ with $S_i$ belong to $B$, and $\pi\cap B=\ell$, necessarily the line $\ell\subseteq S_i$.

For $i \neq j$, the subspace $S_i \cap  S_j$ is either the line $\ell$ or it is a plane that meets $Q$  in two
lines, one of which is the line $\ell$. In particular, at most $2q + 1$ points of $Q$  lie in $S_i \cap  S_j$. It follows
that the union $S_i\cup  S_j$ contains at least $2(q + 1)^2 - (2q + 1) = 2q^2 + 2q + 1$ points of
$B$. Also, as $\ell$ belongs to all the subspaces $S_i$, the union of three of the solids $S_i$ contains at least
$3(q+ 1)^2 - 3(2q + 1) + (q+ 1) = 3q^2 + q+ 1 $ points of $B$.

By (1), we may assume that $e_{P_1} \leq  q- 1$. Consider a point $P_i$ with $i \geq 2$. The line $P_1P_i$ is a
line of $Q$; so it lies in two planes of $Q$. One of these planes is $\pi$. Let $\pi'$ be the other plane.
Then $\pi'$  meets $S_i$ in a line and this line lies in $B$. As $|P_1^\perp \cap B|\leq  (q + 1)^2 + q - 1$, at least
two of the points of this line must lie in $S_1$. Hence, the line $\pi'\cap S_i$ lies in $S_1$. This shows
that $S_1 \cap  S_i$ is a plane that meets $Q$ in two lines $\pi\cap  S_i = \ell$ and $\pi'\cap  S_i$.

The solid $S_1$ lies in $q + 1$ hyperplanes. One of these hyperplanes is $P_1^\perp$. As $|P_1^\perp \cap B| = (q + 1)^2 + e_{P_1} \leq (q+1)^2+q-1<
2q^2+2q+1$ , the hyperplane $P_1^\perp$ cannot contain one of the solids $S_2, \dots, S_{q+2}$. Since $S_1$ spans a
hyperplane with each of the $q+ 1$ solids $S_2,\dots,S_{q+2}$, it follows that there exists a hyperplane $H$
on $S_1$ that contains at least two of the solids $S_2,\dots,S_{q+2}$. Then $|H\cap B|\geq 3q^2 + q + 1$. We
may assume that $S_2,S_3\subseteq H$.

In order to prove the theorem, it suffices to show that $ H \cap B$ meets all the lines of $Q$. Assume
that this is not the case. As $|H \cap B| \geq 3q^2 + q + 1$, \Cref{lemma:block} shows that $H$ is a
non-tangent hyperplane, that $|H \cap B| \leq  4q^2 - 2q+1$, and that every line of $H \cap  Q$  meets
$B$ in at least two points.

Consider the three solids $S_1, S_2, S_3$  of $H$. The line $\ell$ belongs to all three of them. As $H$ is
a non-tangent hyperplane, it meets $Q$  in a quadric of type $Q(4, q)$. Thus $\ell$ meets $ (q + 1)q$
other lines of $H \cap Q=Q(4,q)$. Each of these lines meets $B$ in at least two points. One of these two points
belongs to the line $\ell$. Of the $(q + 1)q$ lines of $H \cap  Q$  that meet $\ell$, at most $3(q + 1)$  lie in one of the solids
$S_1,S_2,S_3$. Therefore, at least $(q + 1)(q - 3)$ lines of $H \cap Q$  meet $\ell$ and are not contained in
$S_1,S_2,S_3$. Since each line of $H$ meets $B$ in at least two points, each of these $ (q+ 1)(q - 3)$ 
lines contains a point of $B$ that is not in $S_1,S_2$ or $S_3$. Different lines yield different points of
$B$, since $H$ does not contain a plane of $Q$ (or, equivalently, $H \cap Q = Q(4,q)$ is a generalised quadrangle) \cite{MR2508121}. Thus we have found $ (q+ 1)(q- 3)$  extra points in
$B$ so that $|H\cap B |\geq 3q^2 + q + 1 + (q + 1)(q- 3)=4q^2-q-2$. But $|B\cap H|\leq  4q^2 - 2q+1$, a contradiction.

\bigskip

{\bf Case 2}. Assume that for all the  points $P\in\pi\setminus \ell$, $e_P=q$.

\bigskip
    
    If $e_P=q$ for every point $P\in\pi\setminus\ell$, then $|P^\perp \cap B|=(q+1)^2+q\leq (q+1)^2+2q-3$. So for every such point $P\in \pi\setminus \ell$, $P^\perp$ contains a solid $S$, with $P\notin S$ and $S\cap Q = Q^+(3,q)\subseteq B$ by Lemma~\ref{lemma:solid}. Moreover, $S\cap\pi=\ell$ for each such solid $S$. We want to find points $P_1,P_2,\dots,P_{q+2}$  such that $S_1\cap S_i$ is a plane for all $i=2,\dots,q+2$, since then we can continue the proof of Case 1.
    
    First choose $q+2$ arbitrary points $P_i\in\pi\setminus\ell$. The line $P_1P_i$ lies in two planes of $Q$, one of which is $\pi$. Let $\pi_i$ be the other plane. Its intersection with $S_i$ is a line $\ell_i$ that is contained in $B$. If this line $\ell_i$ also lies in $S_1$, then $S_1\cap S_i$ is a plane. If this is the case for all $i\in\{2,\dots,q+2\}$, then we are done. So suppose that this is not the case. Then $\ell_i\not\subseteq S_1$ for a certain $i\in\{2,\dots,q+2\}$. W.l.o.g.\ this is the case for $i=2$. Since $\ell_2\subseteq P_1^\perp\cap B$ and $e_P=q$, $P_1^\perp\cap B=Q^+(3,q)_1\cup\ell_2$. Now, for any other point $P_i\in\pi\setminus(\ell\cup\pi_2)$, $\ell_i$ must lie in $S_1$, since otherwise $P_1^\perp\cap B$ contains at least $q-1$ more points (those extra points on the line $\ell_i$) and therefore $e_{P_1}\geq 2q-1$, a contradiction. So $S_1\cap S_i$ is a plane for all points $P_i\in\pi\setminus(\ell\cup\pi_2)$.

    The remainder of the proof is very similar, but replacing \Cref{Res:Metsch:Lemma2.2} by \Cref{lemma:block} and replacing \(4q^2-2q\) by \(4q^2-2q+1\).
\end{proof}

\section{Small weight codewords of \texorpdfstring{$\mc(3,q)$}{C(3,q)}}
 \label{Sec:MinWt}

In this section, we determine the small weight codewords of $\mc(3,q)$.
We do this by investigating the small weight codewords of a code that contains $\mc(3,q)$.

\begin{df}
 Let $C(q)$ denote the code spanned by $\mc(3,q)$ and the characteristic vectors of the lines of the symplectic polar spaces $W(3,q)$ embedded in $\pg(3,q)$.
\end{df}

First, we prove that $C(q) = \mc(3,q)$ if and only if $q$ is even.
We actually prove a similar result for a more general class of codes.

\old{\begin{df}
 For $c: \mg_1(3,q) \to \FF_p$, define
 \[
  \proj(c) : \mg_0(3,q) \to \FF_p:
  P \mapsto \sum_{\substack{\ell \in \mg_1(3,q)\\P\in\ell}} c(\ell).
 \]
\end{df}}
\begin{df}
 Let $0<i<j<n$. For $c: \mg_j(n,q) \to \FF_p$, define
 \[
  \proj^{(i)}(c) : \mg_i(n,q) \to \FF_p:
  \iota \mapsto \sum_{\substack{\lambda \in \mg_j(n,q)\\\iota\subseteq\lambda}} c(\lambda).
 \]
\end{df}

\old{\begin{lm}
A function $c: \mg_1(3,q) \to \FF_p$ belongs to  $\mc(3,q)^\perp$ if and only if $\proj(c)$ belongs to $\mc_{0,1}(3,q)^\perp$.
\end{lm}}

\begin{lm}\label{Lm:Proj}
Let $0<i,j,k<n$, $i<j$. A function $c: \mg_j(n,q) \to \FF_p$ belongs to  $\mc_{j,k}(n,q)^\perp$ if and only if $\proj^{(i)}(c)$ belongs to $\mc_{i,k}(n,q)^\perp$.
\end{lm}

\old{\begin{proof}
 It is valid that $\proj(c)$ belongs to $\mc_{0,1}(3,q)^\perp$ if and only if for each line $\ell$,
 \[
  \sum_{P \in \ell} \proj(c)(P) = 0.
 \]
 On the other hand,
 \begin{align*}
  \sum_{P \in \ell} \proj(c)(P) 
  = \sum_{P \in \ell} \sum_{\substack{\ell' \in \mg_1(3,q) \\ P \in 
 \ell'}} c(\ell')
  = \sum_{\ell' \in \mg_1(3,q)} c(\ell') |\ell \cap \ell'|
  = \sum_{\substack{\ell' \in \mg_1(3,q) \\ \ell \cap \ell' \neq \emptyset}} c(\ell')
  = c \cdot \chi_\ell.
 \end{align*}
Hence, $\proj(c)$ belongs to the dual code of points and lines of $\pg(3,q)$ if and only if $c \cdot \chi_\ell = 0$ for each line $\ell$, which is equivalent to $c \in \mc(3,q)^\perp$.
\end{proof}}
\begin{proof}
 It is valid that $\proj^{(i)}(c)$ belongs to $\mc_{i,k}(n,q)^\perp$ if and only if for each \(k\)-space $\kappa$,
 \[
  \sum_{\substack{\iota\in\mg_i(n,q)\\\iota\cap\kappa\neq\emptyset}} \proj^{(i)}(c)(\iota) = 0.
 \]
 On the other hand,
 \begin{align*}
  \sum_{\substack{\iota\in\mg_i(n,q)\\\iota\cap\kappa\neq\emptyset}} \proj^{(i)}(c)(\iota) 
  = \sum_{\substack{\iota\in\mg_i(n,q)\\\iota\cap\kappa\neq\emptyset}} \sum_{\substack{\lambda\in\mg_j(n,q)\\\iota\subseteq\lambda}} c(\lambda)
  = \sum_{\substack{\lambda\in\mg_j(n,q)\\\lambda\cap\kappa\neq\emptyset}} c(\lambda)
  = c \cdot \chi_\kappa^{(j)},
 \end{align*}
 where we used that
 \[
  | \sett{ \iota \in \mg_i(n,q) }{ \iota \subseteq \lambda, \,  \iota \cap \kappa \neq \emptyset} | \equiv
  \begin{cases}
   0 & \text{if } \lambda \cap \kappa = \emptyset \\
   1 & \text{if } \lambda \cap \kappa \neq \emptyset
  \end{cases}
  \pmod p,
 \]
 because the number of \(i\)-subspaces of a \(j\)-space, having nonempty intersection with a given subspace, equals \(1\) modulo \(p\), see e.g.\ \cite[Lemma 3.2]{bcn}.
 
Hence, $\proj^{(i)}(c)$ belongs to the dual code of \(i\)-spaces and \(k\)-spaces of $\pg(n,q)$ if and only if $c \cdot \chi_\kappa^{(j)} = 0$ for each \(k\)-space $\kappa$, which is equivalent to $c \in \mc_{j,k}(n,q)^\perp$.
\end{proof}

\old{\begin{prop}
 Let $S$ denote the set of absolute lines of a symplectic polar space $W(3,q)$ embedded in $\pg(3,q)$.
 Then $\chi_S \in \mc(3,q)$ if and only if $q$ is even.
\end{prop}}
\begin{prop}
 \label{Prop:Symplectic}
 Let $S$ denote the set of absolute lines of a symplectic polar space $W(2n+1,q)$ embedded in $\pg(2n+1,q)$. Then $\chi_S \in \mc_{1,n}(2n+1,q)$ if and only if $q$ is even.
\end{prop}

\begin{proof}
 \old{First suppose that $q$ is even.
 Consider a hyperbolic quadric $\mq^+(3,q)$, and a regulus $\mathcal R$ of lines contained in $\mq^+(3,q)$.
 Consider the codeword 
 \[
  c = \sum_{\ell \in \mathcal R} \chi_\ell
 \]
 of the code $C(q)$.
 Each point of $\mq^+(3,q)$ is contained in a unique line of $\mathcal R$.
 Take a line $\ell$ in $\pg(3,q)$.
 If $\ell$ intersects $\mq^+(3,q)$ in $k \in \set{0,1,2}$ points, then clearly $c\cdot \chi_\ell = k$.
 Otherwise, $\ell$ is completely contained in $\mq^+(3,q)$.
 Either $\ell$ belongs to $\mathcal R$, and $c\cdot \chi_\ell = 1$, or $\ell$ belongs to the opposite regulus and $c\cdot \chi_\ell = q+1 = 1$.
 Hence, $c$ is the characteristic vector of the singular lines of $\mq^+(3,q)$.
 Since $q$ is even, these are the absolute lines with respect to the polarity corresponding to $\mq^+(3,q)$, which must be symplectic, again by $q$ being even.}

 First suppose that $q$ is even.
 Consider a hyperbolic quadric $\mq^+(2n+1,q)$, and an equivalence class $\mathcal R$ of the generators (\(n\)-spaces) of $\mq^+(2n+1,q)$.
 Consider the codeword 
 \[
  c = \sum_{\lambda\in \mathcal R} \chi^{(1)}_{\lambda}
 \]
 of the code $\mc_{1,n}(2n+1,q)$.
 Take a line $\ell$ in $\pg(2n+1,q)$.
 Then $c(\ell)$ equals the number of generators in $\mr$ that intersect $\ell$ non-trivially.
 Suppose that $\ell$ intersects $Q^+(2n+1,q)$ in $m_1$ points.
 The number of generators of $Q^+(2n+1,q)$ of $\mr$ through a point $P$ equals 0 if $P \notin Q^+(2n+1,q)$ and equals the number of generators of $\mr$ in $P^\perp$ otherwise, with $\perp$ the polarity associated to $Q^+(2n+1,q)$.
 Since $P^\perp$ intersects $Q^+(2n+1,q)$ in a cone $P Q^+(2n-1,q)$, this number equals the number of generators of one equivalence class in $Q^+(2n-1,q)$, which is given by $\prod_{i=1}^{n-1} (q^i+1)$ (where we use the convention that an empty product equals $1$) by \Cref{Res:Quadrics} (\ref{Res:Quadrics:Generators}) and (\ref{Res:Quadrics:Equivalence}).
 Let $m_2$ denote the number of generators of $Q^+(2n+1,q)$ through $\ell$.
 Then
 \[
  c(\ell) = m_1 \prod_{i=1}^{n-1} (q^i + 1) - q m_2 \equiv m_1 \pmod p.
 \]
 Hence, $c$ is the characteristic vector of the singular lines of $Q^+(2n+1,q)$.
 Since $q$ is even, these are the absolute lines with respect to $\perp$, which must be symplectic, again by $q$ being even.

 Now suppose that $q$ is odd.
 Let $\perp$ denote the polarity associated to $W(2n+1,q)$.
 Take a 3-space $\Sigma$ of $W(2n+1,q)$ such that $\Sigma \cap \Sigma^\perp = \emptyset$.
 Then the absolute subspaces of $\perp$ contained in $\Sigma$ are the lines and points of a $W(3,q)$ in $\Sigma$.
 Consider a hyperbolic quadric $Q^+(3,q)$ in $\Sigma$ and let $\mr^+$ and $\mr^-$ denote its two reguli.
 Consider the function
 \[
  c: \mg_1(2n+1,q) \to \FF_p: \ell \mapsto \begin{cases}
   1 & \text{if } \ell \in \mathcal R^+, \\
   -1 & \text{if } \ell \in \mathcal R^-, \\
   0 & \text{otherwise.}
  \end{cases}
 \]
 Then $\proj^{(0)}(c) = \zero$, hence $c \in \mc_{1,n}(2n+1,q)^\perp$ by \Cref{Lm:Proj}.
 Therefore, if $\chi_S$ would be in $\mc_{1,n}(2n+1,q)$, then $\chi_S \cdot c = 0$, which means that $S$ contains the same number of lines of $\mathcal R^+$ and $\mathcal R^-$ modulo $p$.
 Hence, in order to prove that $\chi_S \notin \mc_{1,n}(2n+1,q)$, it suffices to show that there exists a hyperbolic quadric $Q^+(3,q)$ and a symplectic polar space $W(3,q)$ in $\Sigma \cong \pg(3,q)$ where the reguli of $Q^+(3,q)$ do not contain the same number of lines of $W(3,q)$ modulo $p$.
 
 To this end, consider the symplectic polar space of $\pg(3,q)$ defined by the bilinear form
 \[
  b(x,y) = x_0 y_1 - x_1 y_0 + x_2 y_3 - x_3 y_2.
 \]
 From now on, let $\perp$ denote the polarity of the symplectic space $W(3,q)$ corresponding to $b$.
 A line through two distinct points of $\pg(3,q)$, with coordinates $x=(x_0,\dots,x_3)$ and $y=(y_0,\dots,y_3)$, is absolute if and only if $b(x,y) = 0$.
 Now consider the hyperbolic quadric with equation $x_0 x_1 = x_2 x_3$.
 Its two reguli are
 \begin{align*}
  \mathcal R^+ &= \sett{\vspan{(\alpha,0,\beta,0), (0,\beta,0,\alpha)}}{ \vspan{(\alpha,\beta)} \in \pg(1,q)}, \\
  \mathcal R^- &= \sett{\vspan{(\alpha,0,0,\beta), (0,\beta,\alpha,0)}}{\vspan{(\alpha,\beta)} \in \pg(1,q)}.
 \end{align*}
 On the one hand,
 \[
  b\left((\alpha,0,\beta,0),(0,\beta,0,\alpha)\right)
  = 2 \alpha \beta.
 \]
 Since $q$ is odd, this means that $\mathcal R^+$ contains two absolute lines.
 On the other hand
 \[
  b\left((\alpha,0,0,\beta),(0,\beta,\alpha,0)\right) = \alpha \beta - \beta \alpha=0,
 \]
 so all $q+1$ lines in $\mathcal R^-$ are absolute.
 Note that $q+1 \equiv 1 \not \equiv 2 \pmod p$, which finishes the proof.
\end{proof}

Now we continue our investigation of the small weight codewords of $C(q)$.

\begin{lm}
 \label{Lm:Inprod}
 Given $c \in C(q)$, there exists a scalar $\alpha \in \FF_p$ with the property that $c \cdot \chi_S = \alpha$ whenever $S$ is one of the following types of sets:
 \begin{itemize}
     \item $S$ equals the set of all lines in $\pg(3,q)$, i.e.\ $S = \mg_1(3,q)$,
     \item $S$ equals the set of lines in a plane of $\pg(3,q)$,
     \item $S$ equals the set of lines through a point of $\pg(3,q)$,
     \item $S$ equals the set of lines through a point in a plane of $\pg(3,q)$.
 \end{itemize}
\end{lm}

\begin{proof}
 Let $\chi$ either be the characteristic vector of the lines of a symplectic space $W(3,q)$, or a vector of the form $\chi_\ell$ for some line $\ell$ in $\pg(3,q)$.
 Then for each of the sets $S$ listed above, $\chi \cdot \chi_S = 1$.
 The proof for the other codewords of $C(q)$ follows from the  linearity of the code $C(q)$.
\end{proof}

Note that if $\alpha \neq 0$ in the above lemma, then $c \cdot \chi_S \neq 0$ whenever $S$ is the set of lines through a point in a plane.
Therefore, the image of $\supp(c)$ under the Klein correspondence must be a blocking set with respect to the lines of $Q^+(5,q)$.
If $\wt(c)$ is small, then we can apply \Cref{Thm:Main:Klein} to obtain strong structural information about $\supp(c)$.
It remains to deal with the case $\alpha = 0$.

We will use a result in projective planes that is a special case of a result in general combinatorial designs.
For a proof, see \cite[Theorem 3.1]{stinson} for the projective planes and \cite[Theorem 4.1]{stinson} or \cite[Corollary 5.3]{haemers} for the general result.

\begin{res}
 \label{Res:EdgeDom}
 A set $S$ of lines in $\pg(2,q)$ covers at least $\frac{(q+1)^2|S|}{q+|S|}$ points.
\end{res}

We need one more small lemma, and then we are ready to deal with the case $\alpha = 0$.
Given a function $f$, and a subset $S$ of its domain, denote the restriction of $f$ to $S$ by $f_{|S}$.

\begin{lm}
 \label{Lm:Restriction}
 Suppose that $c \in C(q)$.
 \begin{itemize}
  \item If $S$ equals the set of lines in a plane $\pi$ of $\pg(3,q)$, then $c_{|S}$ belongs to the code $\mc_{1,0}(2,q)$ defined on the lines of $\pi$.
  \item If $S$ equals the set of lines through a point $P$, then $c_{|S}$ belongs to the code of points and lines $\mc_{0,1}(2,q)$ in the quotient geometry through $P$.
 \end{itemize}
\end{lm}

\begin{proof}
 Suppose that $S$ denotes the set of lines in the plane $\pi$.
 Then for each line $\ell$, either $\ell \in \pi$ and $\chi^{(1)}_{\ell|S} = \one$, or $\ell \cap \pi$ is some point $P$, and $\chi^{(1)}_{\ell|S}$ is the characteristic vector of the lines through $P$ in $\pi$.
 Similarly, if $\chi$ is the characteristic vector of the lines of a symplectic space $W(3,q)$ with polarity $\perp$, then $\chi_{|S}$ is the characteristic vector of the lines through the point $\pi^\perp$ in $\pi$.
 In all cases, these vectors belong to the code of $\mc_{1,0}(2,q)$ defined on the lines of $\pi$.
 Then this is valid for all $c \in C(q)$ by linearity.

 The proof is analogous for $S$ equal to  the set of lines through a point in $\pg(3,q)$.
\end{proof}

\begin{prop}
 \label{Prop:MinWtPerpToOne}
 Suppose that $q \geq 19$.
 Let $c\in C(q)$, with $c \cdot \one = 0$. If $c \neq \zero$, then $\wt(c)\geq 
 q^3 + 2q^2 + q + 3.$
\end{prop}

\begin{proof}
 Assume that $c \cdot \one = 0$ and $c \neq \zero$.
 Let $\delta$ denote $3$ if $q$ is prime and $12$ otherwise.
 Take a plane $\pi$ which is not disjoint to $\supp(c)$ and let $S$ be the set of lines in $\pi$.
 Consider the code $\mc_{1,0}(2,q)$ of lines and points in $\pi$.
 Then $c_{|S} \in \mc_{1,0}(2,q)$ by \Cref{Lm:Restriction}.
 Now suppose that $\pi$ contains fewer than $3q-\delta$ points of $\supp(c)$.
 By \Cref{Res:PlaneCodes}, $c_{|S}$ is a linear combination of at most two characteristic vectors of points, i.e.\ $c_{|S} = \alpha \chi_{P}^{(1)} + \beta \chi_{R}^{(1)}$ for some scalars $\alpha, \beta \in \FF_p$ and points $P, R \in \pi$.
 \Cref{Lm:Inprod} implies that $c \cdot \chi_S = c \cdot \one = 0$.
 Since $\alpha + \beta = c_{|S} \cdot \one = c \cdot \chi_S = 0$, it follows that $\alpha = -\beta$.
 Note that $\alpha, \beta \neq 0$ and $P \neq R$ since we assumed that $c_{|S} \neq \zero$.
 
 We conclude the following.
 For any plane $\pi$, one of three cases holds:
 \begin{itemize}
  \item $\pi$ contains no lines of $\supp(c)$,
  \item $\pi$ contains exactly $2q$ lines of $\supp(c)$, and there exist two distinct points $P, R$ in $\pi$ so that the lines of $\supp(c)$ contained in $\pi$ are exactly the lines of $\pi$ containing exactly one of the points $P$ and $R$,
  \item $\pi$ contains at least $3q-\delta$ lines of $\supp(c)$.
 \end{itemize}

 Completely analogously, any point of $\pg(3,q)$ belongs to 0, $2q$, or at least $3q-\delta$ lines of $\supp(c)$.

 First consider the case where some plane $\pi$ of $\pg(3,q)$ contains $2q$ lines of $\supp(c)$, and let $P, R$ be as above.
 Then every point of $\pi$ not on the line $\vspan{P,R}$ belongs to two lines of $\supp(c)$ in $\pi$, hence, such a point belongs to at least $2q-2$ lines of $\supp(c)$ outside of $\pi$.
 Also, $P$ and $R$ lie on at least $2q$ lines of $\supp(c)$.
 This yields
 \[
  \wt(c) \geq q^2 (2q-2) + 2\cdot2q =2q^3-2q^2+4q > q^3 + 2q^2 + q + 3.
 \]

 If some point $P$ belongs to exactly $2q$ lines of $\supp(c)$, we can apply a similar argument.
 Thus, the last case to consider is the case where no plane or point of PG$(3,q)$ is incident with exactly $2q$ lines of $\supp(c)$.
 Take a plane $\pi$ containing a line of $\supp(c)$.
 Then $\pi$ contains at least $3q-\delta$ lines of $\supp(c)$.
 By \Cref{Res:EdgeDom}, these lines cover at least $\frac{(q+1)^2(3q-\delta)}{4q-\delta}$ points of $\pi$.
 Each of these points belongs to at least $3q-\delta - (q+1) = 2q - \delta - 1$ lines outside of $\pi$.
 Therefore,
 \[
  \wt(c) \geq 3q-\delta + \frac{(q+1)^2(3q-\delta)}{4q-\delta} (2q-\delta-1).
 \]
 Plugging in that either $\delta = 3$ and $q \geq 19$, or $\delta = 12$ and $q \geq 25$ yields
 \[
  \wt(c) > q^3 + 2q^2 + q + 2. \qedhere 
 \]
\end{proof}

\begin{thm}
 \label{Thm:LargeCode}
 Suppose that $q \geq 19$ and that $c$ is a non-zero codeword of $C(q)$ with $\wt(c) \leq q^3 + 2q^2 + q + 1$.
 Then one of the following cases holds.
 \begin{enumerate}
  \item $\wt(c) = q^3 + q^2 + q + 1$, and $c$ is a scalar multiple of the characteristic vector of the lines of a symplectic space $W(3,q)$.
  \item $\wt(c) = q^3 + 2q^2 + q + 1$, and $c$ is a scalar multiple of $\chi^{(1)}_\ell$ for some line $\ell$.
 \end{enumerate}
\end{thm}

\begin{proof}
 Define $\alpha$ as $c \cdot \one$.
 By \Cref{Prop:MinWtPerpToOne}, $\alpha \neq 0$.
 As noted before, \Cref{Lm:Inprod} implies that $\supp(c)$ is equivalent to a blocking set of $Q^+(5,q)$ under the Klein correspondence, whose size is at most $q^3 + 2q^2 + q + 1$.
 By \Cref{Thm:Main:Klein}, there are two options:
 \begin{itemize}
  \item $\supp(c)$ contains the set $S$ of lines of a symplectic space $W(3,q)$.
  In this case, let $\chi$ denote $\chi_S$.
  Since $S \subseteq \supp(c)$, $\wt(c - \alpha \chi) \leq \wt(c)$.
  \item $\supp(c)$ contains the set $S$ of lines that intersect a given line $\ell$ of $\pg(3,q)$ exactly in a point.
  In this case, let $\chi$ denote $\chi^{(1)}_\ell$.
  Since $\ell$ is the only possible element of $\supp(\chi) \setminus \supp(c)$, $\wt(c - \alpha \chi) \leq \wt(c)+1$.
 \end{itemize}
 In both cases, we find a codeword $c - \alpha \chi$ of $C(q)$ with $(c-\alpha \chi) \cdot \one = 0$ and $\wt(c - \alpha \chi) \leq \wt(c) + 1 < q^3 + 2q^2 + q + 3$.
 By \Cref{Prop:MinWtPerpToOne}, this implies that $c = \alpha \chi$.
\end{proof}

Note that \Cref{Thm:Main:Code} follows directly from \Cref{Thm:LargeCode} and \Cref{Prop:Symplectic}

\paragraph{Acknowledgements.}
To obtain the results of this article, we needed to improve \Cref{Res:Metsch} of Metsch on blocking sets with respect to the lines of the hyperbolic quadric $\mq^+(5,q)$ \cite{metsch2000blocking}. We thank Metsch for allowing us to repeat large parts of his proofs in \Cref{Sec:Klein} of this article to obtain the improvement.
We are also grateful to Mrinmoy Datta for helpful discussions.

Robin Simoens is supported by the Research Foundation Flanders (FWO) through the grant 11PG724N.

\vfill
\noindent\textsc{Sam Adriaensen}\\
\textsc{\small Department of Mathematics and Data Science}\\[-1mm]
\textsc{\small Vrije Universiteit Brussel}\\[-1mm]
\textsc{\small Pleinlaan 2, 1050 Elsene, Belgium}\\
{\it E-mail address:} {\href{mailto:sam.adriaensen@vub.be}{\url{sam.adriaensen@vub.be}}}\\

\noindent\textsc{Robin Simoens}\\
\textsc{\small Department of Mathematics: Analysis, Logic and Discrete Mathematics}\\[-1mm]
\textsc{\small Ghent University}\\[-1mm]
\textsc{\small Krijgslaan 281, 9000 Gent, Belgium}\\
\textsc{\small Department of Mathematics}\\[-1mm]
\textsc{\small Universitat Politècnica de Catalunya}\\[-1mm]
\textsc{\small C. Pau Gargallo 14, 08028 Barcelona, Spain}\\
{\it E-mail address:} {\href{mailto:robin.simoens@ugent.be}{\url{robin.simoens@ugent.be}}}\\

\noindent\textsc{Leo Storme}\\
\textsc{\small Department of Mathematics: Analysis, Logic and Discrete Mathematics}\\[-1mm]
\textsc{\small Ghent University}\\[-1mm]
\textsc{\small Krijgslaan 281, 9000 Gent, Belgium}\\
{\it E-mail address:} {\href{mailto:leo.storme@ugent.be}{\url{leo.storme@ugent.be}}\\

\end{document}